\documentclass[12pt]{article}  
\usepackage{times,amsmath,amsfonts,amsthm,amssymb,eucal}
\usepackage{array}
\usepackage{mathrsfs}
\usepackage{amsmath}
\usepackage{setspace}
\usepackage{fullpage}
\usepackage{fancyhdr}
\usepackage{setspace}
\usepackage{amssymb}
\usepackage{epsfig}
\usepackage{graphicx}
\usepackage{makeidx}
\usepackage{color}
\usepackage{multirow}
\usepackage{epsfig}
\usepackage{graphicx}
\usepackage{makeidx}
\usepackage{gensymb}
\usepackage{enumerate}
\usepackage{multirow}
 \usepackage{setspace}
\usepackage{amssymb}
 \usepackage{setspace}
\newtheorem{theorem}{Theorem}
\newtheorem{proposition}{Proposition}
\newtheorem{lemma}{Lemma}
\newtheorem{corollary}{Corollary}
\def\Theorem{\begin{theorem}\sl}
\def\EndTheorem{\end{theorem}}
\def\Proposition{\begin{proposition}\sl}
\def\EndProposition{\end{proposition}}
\def\Lemma{\begin{lemma}\sl}
\def\EndLemma{\end{lemma}}
\def\Corollary{\begin{corollary}\sl}
\def\EndCorollary{\end{corollary}}

\theoremstyle{remark}

\numberwithin{equation}{section}

\begin{document}
\renewcommand{\baselinestretch}{1.2}
\markright{
\hbox{\footnotesize\rm  2011: Preprint}\hfill
}

\markboth{\hfill{\footnotesize\rm LUAI AL LABADI AND MAHMOUD ZAREPOUR} \hfill}
{\hfill {\footnotesize\rm Dirichlet Process with Large Concentration Parameter} \hfill}

\renewcommand{\thefootnote}{}
$\ $\par


\fontsize{10.95}{14pt plus.8pt minus .6pt}\selectfont
\vspace{0.8pc}
\centerline{\large\bf  THE DIRICHLET PROCESS WITH LARGE}
\vspace{2pt}
\centerline{\large\bf CONCENTRATION PARAMETER}
\vspace{.4cm}
\centerline{LUAI AL LABADI AND MAHMOUD ZAREPOUR}
\vspace{.4cm}

\fontsize{9}{11.5pt plus.8pt minus .6pt}\selectfont

\begin{quotation}
\noindent {\it Abstract:}
 Ferguson's Dirichlet process  plays an important role in  nonparametric Bayesian inference. Let $P_a$ be the Dirichlet process in $\mathbb{R}$ with a base probability  measure $H$ and  a concentration parameter $a>0.$ In this paper, we show that $\sqrt {a} \big(P_a\left((-\infty,t]\right) -H\left((-\infty,t]\right)\big)$ converges to a certain Brownian bridge as $a \to \infty.$ We also derive a certain  Glivenko-Cantelli theorem for the Dirichlet process. Using the functional delta method, the weak convergence of the quantile process is also obtained.
A large concentration parameter occurs when a statistician puts too much emphasize on his/her prior guess. This scenario also happens  when the sample size is large and the posterior is used to make inference.\par

\vspace{9pt}

\noindent {\it Key words and phrases:}
Bayesian nonparametric, Brownian bridge, Dirichlet process, quantile process, weak convergence.
\par

\end{quotation}\par


\fontsize{10.95}{14pt plus.8pt minus .6pt}\selectfont
\pagestyle{myheadings}\markboth{}{LUAI AL LABADI AND MAHMOUD ZAREPOUR}

\section{Introduction}
\label{intro}

 In nonparametric Bayesian inference, we need to place a prior on an infinite dimensional space such as the space of probability measures. Ferguson (1973) used a Dirichlet process (a normalized gamma process) as a prior on this space. For $k \ge 2,$ we say that the random vector  $(Y_1,\ldots,Y_{k})$ has the Dirichlet distribution
with parameters $(a_1,\ldots,a_{k})$, where $a_i > 0$ for all $i,$ if it has the joint probability density function
\begin{equation}
f(y_1,\ldots,y_{k})=\frac {\Gamma\left(\sum_{i=1}^k a_i\right)}{\prod_{i=1}^k \Gamma\left(a_i\right)}\prod_{i=1}^{k} y_i^{a_i-1}~
I_\mathbb{S}(y_1,\ldots,y_{k}) , \nonumber
\end{equation}
where $\mathbb{S}=\left\{ (y_1,\ldots,y_{k}): y_i \ge 0, \sum_{i=1}^{k}y_i=1\right\}$ and $\Gamma(x)=\int_{0}^\infty t^{x-1}e^{-t}dt,\ \ x>0.$ We denote by $D(a_1,\ldots,a_k)$ the Dirichlet distribution with parameters $a_1,\ldots,a_k.$\\

The Dirichlet process was defined in Ferguson (1973) as follows: let $(\mathfrak{X},\mathcal{A})$ be an arbitrary
measurable space and $H$ be a probability measure on  $(\mathfrak{X},\mathcal{A}).$ Let $a>0$ be arbitrary. A random probability
measure $P_a=\left\{P_a(A)\right\}_{A \in \mathcal{A}}$ is called a
Dirichlet process on $(\mathfrak{X},\mathcal{A})$ with parameters
$a$ and $H$, if for any finite measurable partition $\{A_1, \ldots,
A_k\}$ of $\mathfrak{X}$, the joint distribution of the vector
$\left(P_a(A_1), \ldots\,P_a(A_{k})\right)$ has the Dirichlet distribution
with parameters $(a H(A_1), \ldots,$ $aH(A_k)).$ The subscript $a$ is added since in the forthcoming sections we will study  the asymptotic behavior of the random probability measure $P_a$ for large values of  $a$.  We assume that if $H(A_j)=0$, then $P_a(A_j)=0$
with probability one. We write $P\sim \text{DP}(a, H)$ to denote the Dirichlet process with parameters $a$ and $H.$ Throughout this paper, we use the same letter for the probability measure and  its corresponding cumulative distribution function, i.e. $P_a(t)=P_a\left((-\infty,t]\right)$ and $H(t)=H\left((-\infty,t]\right)$. We also assume that the cumulative distribution function $H$ is continuous.\\

For any $A \in \mathcal{A},$ $P_a(A)$ has a Beta distribution
with parameters $a H(A)$ and $a(1-H(A))$. Thus,
\begin{equation}
{E}(P_a(A))=H(A) \ \ \text{ and }\ \ {Var}(P_a(A))=\frac{H(A)(1-H(A))}{1+a}. \label{As.10}
\end{equation}
Furthermore,  for any two sets $A_i$ and $A_j \in \mathcal{A},$ it follows from the properties of a Dirichlet distribution that (Wilks 1963, page 177)
\begin{equation}
{E}(P_a(A_i)P_a(A_j))=\frac{a}{1+a}H(A_i)H(A_j) \label{asy20}
\end{equation}

The probability measure $H$ is called the {base measure} of $P_a.$ Clearly, form (\ref{As.10}), $H$ plays the role of the center
 of the process, while  $a$  can be viewed as  the concentration parameter. The larger $a$ is, the more likely it is that the realization of $P$ is close to $H$. Specifically, for any fixed set  $A \in \mathcal{A}$ and $\epsilon>0,$ we have $P_a(A) \overset{p} \to H(A)$ as $a \to \infty$ since
\begin{equation}
\Pr\left\{\left|P_a(A)-H(A)\right|>\epsilon\right\} \le \frac{H(A)(1-H(A))}{\epsilon^2(1+a)}.\nonumber
\end{equation}
In this paper, ``$\overset{p} \to$"  denotes the convergence in probability.\\

An attractive property of the Dirichlet process is its conjugacy property. That is, if  $X_{1},\ldots, X_{n}$ is a random sample from $P_a\sim DP(a, H)$, then the posterior distribution of $P_a$ given $ X_{1},\ldots, X_{n}$ coincides with the distribution of the Dirichlet process with parameter measure $a^{*} H^{*}$, where
\begin{equation}
a^{*}=a+n \quad \text{and} \quad H^{*}=\frac{a }{a+n}H+\frac{n}{a+n}\frac{\sum_{k=1}^{n}\delta_{{X}_k}}{n}.\label{eq71}
\end{equation}
Here and throughout the paper $\delta_X$ denotes the Dirac measure at $X$, i.e. $\delta_X(A)=1$ if $X \in A$ and $0$ otherwise.\\

Notice that the posterior base distribution $H^{*}$ is a convex combination of the base distribution
 and the empirical distribution. The weight associated with the prior base distribution $H$ is
 proportional to $a$, giving another reason to call $a$ the concentration parameter. The weight
 associated with the empirical distribution is proportional to the number of observations $n$.
The posterior base distribution $H^*$  approaches the prior base measure $H$ for large values of $a.$ On the
 other hand, for small values of $a $, $H^*$ is close to the empirical process.\\

The Dirichlet process has the following  series representation:
\begin{equation}
P_a(\cdot)=\sum_{i=1}^\infty J_i \delta_{\theta_i} (\cdot), \label{asymp:01}
\end{equation}
where $(\theta_i)_{i \geq 1}$  is a sequence of independent and identically distributed (i.i.d.) random variables with common distribution $H$ and $(J_i)_{i \geq 1}$ are random variables chosen to be independent of $(\theta_i)_{i \geq 1}$ and such $0 \le J_i \le 1$ and $\sum_{i=1}^\infty J_i=1$ almost surely.  For several representations for $(J_i)_{i \geq 1},$  see, for example, Ferguson, Phadia, and Tiwari (1992). It follows from (\ref{asymp:01}) that any realization of the Dirichlet process must be a discrete probability measure.\\

Sethuraman and Tiwari (1982) studied the convergence and tightness of Dirichlet processes as the parameters are allowed to converge in  a certain sense. They showed that as the concentration parameter $a \to 0,$ the Dirichlet process converges to a degenerate probability measure at a particular point in  $\mathfrak{X}$ randomly chosen from $H.$  \\

Let $\mathscr{S}$  be a collection of Borel sets in $\mathbb{R}.$ For large values of the concentration parameter $a,$ we study  the weak convergence of the centralized and scaled Dirichlet process defined by
\begin{equation}
D_a(S)=\sqrt{a}\left(P_a(S)-H(S)\right), \ \ S \in \mathscr{S}. \label{eq2.18}
\end{equation}

We also derive the limiting distribution of the Dirichlet quantile process
\begin{equation}
Q_a=\sqrt{a}\left(P^{-1}_a-H^{-1}\right),\label{Asy:eq102}
\end{equation}
where in general the inverse of a distribution function $F$ is given by
$F^{-1}(t)=\inf\left\{x:F(x)\ge t \right\}.$ Moreover, a certain Glivenko-Cantelli theorem for the Dirichlet process for large values of concentration parameter is obtained.\\

For the Dirichlet posterior processes with parameters given in (\ref{eq71}), the concentration parameter $a^* \to \infty$ whenever $n \to \infty$ ($n$ is the sample size). Lo (1987) studied completely the behavior of the process
$$d_{n,a}(t)=\sqrt{n}\left(P^*_{n,a}(t)-F_n(t)\right), \ \ t \in \mathbb{R},$$
as the sample size $n$ gets large, where $P^*_{n,a}$ is the posterior of the Dirichlet process $P_a$ given the data and $F_n$ is the empirical distribution function.  Using this result, Lo (1987) gave an asymptotic justification of the use of Bayesian bootstrap and provided large sample  Bayesian bootstrap probability intervals for the mean, the variance, and bands for the distributions.

\section{Asymptotic Properties of the Dirichlet process}
In this section,  we study the asymptotic properties of  $P_a$ as $a \to \infty,$ where  $P_a\sim DP(a,H).$ Since $H$ is strictly increasing, we have $$\theta_i < t \ \ \text{if and only if} \ \ H(\theta_i) < H(t).$$
Thus,
$$P_a\left((-\infty,t]\right)=\sum_{i=1}^{\infty} J_i \delta_{\theta_i}\left((-\infty,t]\right)\overset{d}=\sum_{i=1}^{\infty} J_i \delta_{H(\theta_i)}\left((0,H(t)]\right).$$
Throughout this paper, ``$\overset{d}=$"  denotes equality in distribution. Since $(\theta_i)_{i \geq 1}$ is a sequence of
i.i.d. random variables with continuous distribution $H,$ for $i \ge 1,$ $U_i\overset{d}=H(\theta_i)$ where $\left\{U_i\right\}_{i \geq 1}$is a sequence of i.i.d. random variables with a uniform distribution on $[0,1]$. Hence,
$$P_a\left((-\infty,t]\right)\overset{d}=\sum_{i=1}^{\infty} J_i \delta_{U_i}\left((0,H(t)]\right).$$
Therefore, without loss of generality, we only consider the case when $H(t)=t$ (i.e., $(\theta_i)_{i \geq 1}$ is a sequence of i.i.d. random variables with uniform distribution on $[0,1]$). Hence,  the process in (\ref{eq2.18}) reduces to
\begin{eqnarray}
\label{Asy:100} D_a(S)=\sqrt{a}\left(P_a(S)-\lambda(S)\right),
 \end{eqnarray}
where $\lambda$ is the Lebesgue measure  on $[0,1].$ Hereafter, unless otherwise stated, $P_a\sim DP(a,\lambda),$ where $\lambda$ is the Lebesgue measure  on $\mathfrak{X}=[0,1].$\\

We now recall the definition of a Brownian bridge indexed by $\mathscr{S}.$ A Gaussian process $\left\{B_\lambda(S): S \in \mathscr{S}\right\}$ is called a Brownian bridge if
\begin{eqnarray}
\label{asym:cov0} {E}\left[B_\lambda(S)\right]=0 \quad \text {and} \quad {Cov}\left(B_\lambda(S_i),B_\lambda(S_j)\right)=\lambda(S_i \cap S_i)- \lambda(S_i)\lambda(S_i),
 \end{eqnarray}
where $S, S_i, S_j \in \mathscr{S}$ (Massart 1989).\\

\vspace{2mm}
The next lemma gives the limiting distribution of the process  (\ref{Asy:100}) for any finite Borel sets $S_1,\ldots ,S_k \in \mathscr{S}.$ The proof of the lemma for $k=2$ is given in the appendix and it can be generalized easily to the case of arbitrary $k.$ In this paper, ``$\overset{d} \to$"  denotes the convergence in distribution.

\Lemma \label{Asy.L2} Let $D_a$ be as defined in  (\ref{Asy:100}). Then, as $a \to \infty,$ for any fixed  sets $S_1,\ldots ,S_k \in \mathscr{S}$,
\begin{equation}
\left(D_a(S_1),D_a(S_2),\ldots,D_a(S_k)\right)\overset {d }\to \left(B_\lambda(S_1),B_\lambda(S_2),\ldots,B_\lambda(S_k)\right), \nonumber
\end{equation}
where  $B_\lambda$ is the Brownian bridge indexed by $\mathscr{S}$  with the mean and the covariance structure as given in  (\ref{asym:cov0}).
\EndLemma

\noindent {\textbf{Remark 1.}} The convergence obtained in Lemma \ref{Asy.L2} is called convergence in total variation. This type of convergence is stronger than convergence in distribution (Billingsley 1999, page 29).\\

\noindent {\textbf{Remark 2.}} It follows from Lemma \ref{Asy.L2} that, for any fixed Borel set $S \in \mathscr{S},$
\begin{equation}
D_a(S)=\sqrt {a} \left(P_a(S)-\lambda(S)\right) \overset{d} \to B_\lambda(S),\nonumber
\end{equation}
where $B_\lambda(S)$ is distributed as $N(0,\lambda(S)(1-\lambda(S))).$

\vspace{3mm}
 Lemma \ref{Asy.L2} proves that the finite-dimensional  distributions of the process $D_a$ converge to the corresponding finite-dimensional distribution of $B_\lambda.$   The next theorem shows that the process $D_a$ converges to  the process $B_\lambda$ on    $D[0,1]$ with respect to the Skorokhod  topology.

\Theorem \label{Asy.T1} Let $D_a$ be as defined in  (\ref{Asy:100}). Then , as $a \to  \infty,$ we have:
 \begin{eqnarray}
 \nonumber \sqrt{a}\left(P_a(\cdot)-\lambda(\cdot)\right) \overset{d} \to B_\lambda(\cdot)
 \end{eqnarray}
 on $D[0,1]$ with respect to the Skorokhod  topology, where $B_\lambda$ is a Brownian bridge.
\EndTheorem

\proof
From Lemma \ref{Asy.L2} and  Theorem 13.5 of Billingsley (1999) we need only to prove that for any $0 \le t_1\le t \le t_2 \le 1,$
 \begin{eqnarray}
 \nonumber E\left[\left|P_a(t)-P_a(t_1)\right|^{2\beta}\left|P_a(t_2)-P_a(t)\right|^{2\beta}\right] \le \left|F(t_2)-F(t_1)\right|^{2\alpha},
 \end{eqnarray}

for some  $\beta \ge 0,$ $\alpha>1/2,$ and a nondecreasing continuous function $F$ on $[0,1].$  Take $\beta=1/2,$ $\alpha=1,$ and $F(t)=t$ to show that:
\begin{eqnarray}
 \label{As.100}       E\left[\left(P_a(t)-P_a(t_1)\right)\left(P_a(t_2)-P_a(t)\right)\right] \le \frac{a}{a+1}\left(t_2-t_1\right)^{2}.
 \end{eqnarray}
 Observe that,
 \begin{eqnarray}
 \nonumber \left(P_a(t)-P_a(t_1),P_a(t_2)-P_a(t)\right)\sim D\left(a\lambda(t_1,t],a\lambda(t,t_2],a\left(1-\lambda(t_1,t]-\lambda(t,t_2]\right)\right)
 \end{eqnarray}
 From (\ref{asy20}) we have:
 \begin{eqnarray}
 \nonumber E\left[\left(P_a(t)-P_a(t_1)\right)\left(P_a(t_2)-P_a(t)\right)\right] &=&\frac{a}{a+1}\lambda(t_1,t]\lambda(t,t_2]\\
\nonumber                                                                  &=&\frac{a}{a+1}\left(t-t_1\right)\left(t_2-t\right).
 \end{eqnarray}
Since  $t_1\le t \le t_2,$   (\ref{As.100}) follows. This completes the proof of the theorem.
\endproof

As in Ferguson (1973), under the squared error loss and Dirichlet prior, the no data estimate (or the posterior estimate) for the distribution is the prior distribution $H$. Under the absolute deviation loss, the estimate will be the median of the Dirichlet process with the prior distribution of $H$. Therefore, the Dirichlet quantile process plays a role in estimation. The following corollary derives the asymptotic behavior of the Dirichlet quantile process defined by  (\ref{Asy:eq102}) when the concentration parameter $a$ is large.

\Corollary \label{Asy.} Let $0<p<q<1, $ and  $H$ be a continuous function with positive derivative $h$ on the interval $\left[H^{-1}(p)-\epsilon,H^{-1}(q)+\epsilon\right]$ for some $\epsilon>0.$ Let $Q_a$ be the Dirichlet quantile process defined in  (\ref{Asy:eq102}), where $P_a\sim DP(a,H).$ Then, as $a \to  \infty,$ we have:
 \begin{eqnarray}
 \nonumber Q_a(\cdot) \overset{d} \to -\frac{B_\lambda(\cdot)}{h(H^{-1}(\cdot))}=Q(\cdot),
 \end{eqnarray}
 in $D[p,q].$ That is, the limiting process is a Gaussian process with zero-mean and covariance function
 $$Cov\left(Q(S_i),Q(S_j)\right)=\frac{\lambda(S_i \cap S_j)-\lambda(S_i)\lambda(S_j)}{h(H^{-1}(S_i))h(H^{-1}(S_j))},\quad S_i, S_j \in \mathscr{S}.$$
\EndCorollary

\proof
By Theorem \ref{Asy.T1}  the process $\sqrt{a}\left(P_a-H\right)$ converges in distribution to the process $B_H=B_\lambda(H)=B_\lambda \circ H.$ Almost all sample paths of the limiting process are continuous  on the interval $\left[H^{-1}(p)-\epsilon,H^{-1}(q)+\epsilon\right].$  By Lemma 3.9.23. page 386  of Van der Vaart and Wellner (1996),  the inverse map $H \mapsto H^{-1}$ is Hadamard-differentiable at $H$ tangentially to the subspace of functions that are continuous on this interval. By the functional delta method (Theorem 3.9.4 page 374 of Van der Vaart and Wellner (1996)) we have
\begin{eqnarray}
 \nonumber Q_a(\cdot) \overset{d} \to -\frac{B_\lambda\circ H\circ H^{-1}(\cdot)}{h(H^{-1}(\cdot))}=-\frac{B_\lambda(\cdot)}{h(H^{-1}(\cdot))}
 \end{eqnarray}
 in $D[p,q]$. This completes the proof of the corollary.
\endproof

 \noindent {\textbf{Remark 3.}}  Paralleling Remark 1 of Bickel and Freedman (1981), if  $H^{-1}(0+)>-\infty$ and $H^{-1}(1)<\infty$ and $h$ is continuous on $[H^{-1}(0+),$ $H^{-1}(1)],$ the conclusion of the corollary holds in  $D\left[H^{-1}(0+),H^{-1}(1)\right].$  For example, if $H$ is a uniform distribution on $[0,1],$ then the convergence holds in $D[0,1].$ More generally, we may have one end of the support finite and the other infinite and a modified form of Corollary 1 still  holds. Also from the result of Theorem \ref{Asy.T1} , we can derive asymptotic properties of any  Hadamard-differentiable functional of the $DP(a,H)$ as $a\to \infty.$\\

\noindent {\textbf{Example 1 (Median). }}   Let $M_a$ be the median of $P_a$ and $m$ be the median of $H$ (i.e. $P^{-1}_{a}(0.5)=M_{a}$ and $H^{-1}(0.5)=m$). From Corollary \ref{Asy.} we have:
\begin{eqnarray}
 \nonumber \sqrt{a}\left(M_a-m\right) \overset{d} \to N\left(0,\frac{1}{4h^2(m)}\right),
 \end{eqnarray}
 where $h=H^{\prime}.$ Note that, the asymptotic distribution  of the median for Dirichlet process coincide with that of  the sample median.\\

\noindent {\textbf{Example 2 (Interquantile Range).}}  Similar to Example 1, let  $IQR=Q_{3,a}-Q_{1,a},$ where  $Q_{3,a}$ and $Q_{1,a}$ are  the third and the first quartiles of $P_a$ (i.e. $P^{-1}_{a}(0.75)=Q_{3,a}$ and $P^{-1}_{a}(0.25)=Q_{1,a}$). Let $q_{3}$ and $q_{1}$ be  the third and the first quartiles of $H.$ From Corollary \ref{Asy.},  a simple calculation shows
 \begin{eqnarray}
 \nonumber \sqrt{a}\left(IQR-(q_3-q_1)\right) \overset{d} \to N\left(0,\frac{3}{h^2(q_3)}+\frac{3}{16h^2(q_1)}-\frac{2}{h(q_1)h(q_3)}\right).
 \end{eqnarray}
This gives with the asymptotic distribution of the sample interquartile range.\\

In the  next theorem we establish the Glivenko-Cantelli theorem for the Dirichlet process. In this paper, ``$\overset{a.s.} \to$"  denotes  the almost sure convergence.

\Theorem \label{Asy.T2} Let $P_a\sim DP(a,H)$. Then,
 \begin{eqnarray}
 \nonumber \sup_{x \in \mathbb{R}}\left|P_a(x)-H(x)\right| \overset{a.s.} \to 0,
 \end{eqnarray}
 as $a \to  \infty.$
\EndTheorem

 \proof
From Donoho and Liu (1988),
\begin{eqnarray}
 \label{Asy:DCT} \frac{\left(\sup_x\left|P_a(x)-H(x)\right|\right)^{3/2}}{3^{1/2}} \le  \int_{-\infty}^{\infty}\left(P_a(x)-H(x)\right)^2dH(x).
 \end{eqnarray}
 Notice that $P_a(x) \overset{p} \to H(x),$ as $a \to \infty,$ and $\left(P_a(x)-H(x)\right)^2$ is dominated by 1. Thus, by the dominated convergence theorem  (which remains valid for convergence in probability (Royden 1968, page 92)), we obtain that the right hand side of ( \ref{Asy:DCT}) converges to zero.
  \endproof

When the concentration parameter is large, the Dirichlet process and its corresponding quantile process   share many asymptotic properties with the empirical process and the quantile process.

\section{Acknowledgments.} The authors would like to thank Professor Raluca Balan for her helpful comments and suggestions. The research of the authors are supported by research grants from the \textbf{Natural Sciences and
Engineering Research Council of Canada (NSERC)}.

\vskip .65cm
\noindent
 Luai Al-Labadi is a PhD Student, Department of Mathematics and
Statistics, University of Ottawa, Ottawa, Ontario, K1N 6N5, Canada,
\vskip 2pt
\noindent
E-mail: (lalla046@uottawa.ca)
\vskip 2pt

\noindent
Mahmoud Zarepour is Associate Professor, Department of Mathematics and
Statistics, University of Ottawa, Ottawa, Ontario, K1N 6N5, Canada,
\vskip 2pt
\noindent
E-mail: (zarepour@uottawa.ca)
\vskip .3cm

\newpage
 \section*{Appendix: Proof of Lemma 1 for $\mathbf{k=2}$ }

Assume that $S_1 \cap S_2 =\emptyset.$ (The general case when $S_1$ and $S_2$ are not necessarily disjoint follows from the continuous mapping theorem). Note that
 \begin{eqnarray}
\nonumber \left(P_a(S_1),P_a(S_2),1-P_a(S_1)-P_a(S_2)\right) &\sim &D\big(a\lambda(S_1),a\lambda(S_2),\\ \nonumber  &&a(1-\lambda(S_1)-\lambda(S_2))\big)
 \end{eqnarray}
Set  $X_{i,a}=P_a(S_i)$ and $l_i=\lambda(S_i),~ i=1,2.$ Thus, the joint density function of $P_{1,a}$ and $P_{2,a}$ is:
\begin{equation}
f_{X_{1,a},X_{2,a}}(x_1,x_2)=\frac{\Gamma(a)} {\Gamma(al_1)\Gamma(al_2)\Gamma(a(1-l_1-l_2))} x_1^{al_1-1}x_2^{al_2-1}(1-x_1-x_2)^{a(1-l_1-l_2)-1}.\nonumber
\end{equation}
The joint probability density function of $D_{1,a}=\sqrt {a} \left(X_{1,a}-l_1\right)=\sqrt{a}\left(P_a(S_1)-\lambda(S_1)\right)$ and  $D_{2,a}=\sqrt {a} \left(X_{2,a}-l_2\right)=\sqrt{a}\left(P_a(S_2)-\lambda(S_2)\right)$ is:
\begin{eqnarray}
\nonumber f_{D_{1,a},D_{2,a}}(y_1,y_2)&=&\frac{\Gamma(a)} {a\Gamma(al_1)\Gamma(al_2)\Gamma(a(1-l_1-l_2))} \left(\frac{1}{\sqrt{a}}y_1+l_1\right)^{al_1-1} \\ \nonumber && \left(\frac{1}{\sqrt{a}}y_2+l_2\right)^{al_2-1}\left(1-\frac{y_1+y_2}{\sqrt{a}}-l_1-l_2\right)^{a(1-l_1-l_2)-1}.
\end{eqnarray}
By Scheff\'e's theorem (Billingsely 1999, page 29), it is enough to show that:
\begin{eqnarray}
\label{asy4} f_{D_{1,a},D_{2,a}}(y_1,y_2) \to f(y_1,y_2) =\frac{1}{2\pi |\Sigma|^{1/2}}\exp\left\{-(y_1~y_2)\Sigma^{-1}(y_1~y_2)^T/2\right\},
\end{eqnarray}
where $\Sigma=\begin{bmatrix}  l_1\left(1-l_1\right)& -l_1l_2 \\ -l_1l_2 & l_2\left(1-l_2\right) \end{bmatrix}.$\\
\vspace{1mm}

Use Stirling's formula (Wilks 1963, page 177)
\begin{eqnarray}
\nonumber \Gamma(z)\approx\sqrt{2\pi}z^{z-\frac{1}{2}}e^{-z}, \text{ as } z\to \infty,
 \end{eqnarray}
 where we use the notation $f(z) \approx g(z)$ as $z \to \infty$  if $\lim_{z \to \infty} \frac{f(z)}{g(z)}=1,$
to get:
\begin{eqnarray}
\nonumber \lim_{a \to \infty}f_{D_{1,a},D_{2,a}}(y_1,y_2)&=& \frac {1}{{2\pi}}\lim_{a \to \infty} \Biggl[\frac{\left(\frac{1}{\sqrt{a}}y_1+l_1\right)^{al_1-1}\left(\frac{1}{\sqrt{a}}y_2+l_2\right)^{al_2-1}} {l_1^{al_1-\frac{1}{2}} l_2^{al_2-\frac{1}{2}}}\\
\nonumber &&\frac {\left(1-\frac{1}{\sqrt{a}}y_1-\frac{1}{\sqrt{a}}y_2-l_1-l_2\right)^{a(1-l_1-l_2)-1}} {(1-l_1-l_2)^{(1-l_1-l_2)a-\frac{1}{2}}}\Biggl]\\
\nonumber &=& \frac {1}{2\pi \sqrt{l_1 l_2(1-l_1-l_2)}}\lim_{a \to \infty} \Biggl[\frac{\left(\frac{1}{\sqrt{a}}y_1+l_1\right)^{al_1-1}} {l_1^{al_1-1}}\\
\nonumber && \frac{\left(\frac{1}{\sqrt{a}}y_2+l_2\right)^{al_2-1}\left(1-\frac{1}{\sqrt{a}}y_1-\frac{1}{\sqrt{a}}y_2-l_1-l_2\right)^{a(1-l_1-l_2)-1}} {l_2^{al_2-1}(1-l_1-l_2)^{a(1-l_1-l_2)-1}}\Biggl]\\
\nonumber &=& \frac {1}{2\pi \sqrt{l_1 l_2(1-l_1-l_2)}}\lim_{a \to \infty} \Biggl[\left(1+\frac{y_1}{\sqrt{a}l_1}\right)^{al_1}\\
\nonumber &&\left(1+\frac{y_2}{\sqrt{a}l_2}\right)^{al_2}\left(1-\frac{y_1+y_2}{\sqrt{a}(1-l_1-l_2)}\right)^{a(1-l_1-l_2)}\Biggl]\\
\label{asy6} &=& \frac {1}{2\pi \sqrt{\sigma_{11}\sigma_{22}(1-\rho_{12}^2)}}\exp\left\{\lim_{a \to \infty} a \ln  v_a\right\},
 \end{eqnarray}
where \begin{eqnarray}
\label{asy5} \sigma_{11}=l_1(1-l_1),\ \ \sigma_{22}=l_2(1-l_2), \ \  \rho_{12}=-\sqrt{\frac{l_1l_2}{(1-l_1)(1-l_2)}},
\end{eqnarray}
and
$$v_a=\left(1+\frac{y_1}{\sqrt{a}l_1}\right)^{l_1}\left(1+\frac{y_2}{\sqrt{a}l_2}\right)^{l_2}\left(1-\frac{y_1+y_2}{\sqrt{a}(1-l_1-l_2)}\right)^{1-l_1-l_2}.$$

Observe that,
\begin{eqnarray}
\nonumber\lim_{a \to \infty} a\ln v_a &=&\lim_{a \to \infty}  \frac {1}{1/a} \Biggl[l_1\ln\left(1+\frac{y_1}{\sqrt{a}l_1}\right)+l_2\ln\left(1+\frac{y_2}{\sqrt{a}l_2}\right)\\
\nonumber &&  +(1-l_1-l_2) \ln \left(1-\frac{y_1+y_2}{\sqrt{a}(1-l_1-l_2)}\right)\Biggl]
\end{eqnarray}
 Using L'Hospital's rule we obtain $\lim_{a \to \infty} a\ln v_a$ equals to:
 \begin{eqnarray}
\nonumber && \lim_{a \to \infty} \left[\frac{l_1\frac {-y_1}{2l_1 a^{3/2}}}{\left(1+\frac{y_1}{l_1\sqrt{a}}\right)}+\frac{l_2\frac {-y_2}{2l_2 a^{3/2}}}{\left(1+\frac{y_2}{l_2\sqrt{a}}\right)}-\frac{(1-l_1-l_2) \frac {-(y_1+y_2)}{2(1-l_1-l_2)a^{3/2}}}{\left(1-\frac{y_1+y_2}{(1-l_1-l_2)\sqrt{a}}\right)}\right] (-a^2)\\
\nonumber &=& \lim_{a \to \infty} \frac{a}{2}\left[\frac{l_1y_1}{l_1\sqrt{a}+y_1}+\frac{l_2y_2}{l_2\sqrt{a}+y_2}-\frac{(1-l_1-l_2)(y_1+y_2)}{(1-l_1-l_2)\sqrt{a}-(y_1+y_2)}\right]\\
\nonumber &=&  \lim_{a \to \infty} \frac{a}{2}\Biggl[\frac{l_1y_1y_2+(1-l_2)y_1^2}{\left(l_1\sqrt{a}+y_1\right)\left((1-l_1-l_2)\sqrt{a}-(y_1+y_2)\right)}\\
\nonumber && + \frac{l_2y_1y_2+(1-l_1)y_2^2}{\left(l_2\sqrt{a}+y_2\right)\left((1-l_1-l_2)\sqrt{a}-(y_1+y_2)\right)}\Biggl]\\
\nonumber &=& -\frac{l_2(1-l_2)y_1^2+2l_1l_2y_1y_2+l_1(1-l_1)y_2^2}{2l_1l_2(1-l_1-l_2)}\\
\nonumber &=&- \frac {(1-l_1)(1-l_2)}{2(1-l_1-l_2)}\Biggl[\left(\frac{y_1}{\sqrt{l_1(1-l_1)}}\right)^2+\left(\frac{y_2}{\sqrt{l_2(1-l_2)}}\right)^2\\
\nonumber && +\frac{2y_1y_2}{{(1-l_1)(1-l_2)}}\Biggl]\\
\nonumber &=&-\frac{1}{2(1-\rho_{12}^2)}\Biggl[\left(\frac{y_1}{\sqrt{\sigma_{11}}}\right)^2+\left(\frac{y_2}{\sqrt{\sigma_{22}}}\right)^2-2\rho_{12} \left(\frac{y_1}{\sqrt{\sigma_{11}}}\right) \left(\frac{y_1}{\sqrt{\sigma_{11}}}\right)\Biggl],
\end{eqnarray}
where
$\sigma_{11}, \sigma_{22}$ and $\rho_{12}$ are defined in (\ref{asy5}). The proof follows by using (\ref{asy6}).

\end{document}